\documentclass[ijoc]{informs4anon}

\usepackage{booktabs}
\usepackage{fancyvrb}
\usepackage{multirow}
\usepackage{pgfplots}
\usepackage{standalone}
\usepackage{tikz}
\usepackage{url}

\pgfplotsset{
    compat = 1.8
}
\usepgfplotslibrary{groupplots}

\definecolor{revision_color}{rgb}{0.0, 0, 0}
\newcommand{\revision}{}

\definecolor{revision_color_b}{rgb}{0.0, 0, 0}
\newcommand{\revisionb}{}

\usetikzlibrary{
    arrows,
    shapes,
}

\tikzset{
    ->,
    my_node/.style={
        minimum height=1.5cm,
        minimum width=3cm,
        draw,
    },
    mpb_node/.style={
        minimum height=3cm,
        minimum width=12cm,
        draw,
    },
}

\usepackage{natbib}
 \bibpunct[, ]{(}{)}{,}{a}{}{,}
 \def\bibfont{\small}

\TheoremsNumberedThrough
\EquationsNumberedThrough

\begin{document}


\RUNAUTHOR{Legat, Dowson, Garcia, and Lubin}
\RUNTITLE{MathOptInterface: a data structure for mathematical optimization problems}

\TITLE{MathOptInterface: a data structure for mathematical optimization problems}

\ARTICLEAUTHORS{%
\AUTHOR{Beno\^it Legat}
\AFF{ICTEAM, UCLouvain, Louvain-la-Neuve, Belgium, \EMAIL{benoit.legat@uclouvain.be}}
\AUTHOR{Oscar Dowson}
\AFF{Department of Industrial Engineering and Management Sciences at Northwestern University, Evanston, IL, \EMAIL{o.dowson@gmail.com}}
\AUTHOR{Joaquim Dias Garcia}
\AFF{PSR \& PUC-Rio, Rio de Janeiro, Brazil, \EMAIL{joaquimgarcia@psr-inc.com}}
\AUTHOR{Miles Lubin}
\AFF{Google Research, New York, NY, \EMAIL{mlubin@google.com}}
}

\ABSTRACT{
\revision{We introduce MathOptInterface, an abstract data structure for representing mathematical optimization problems based on combining pre-defined functions and sets. MathOptInterface is significantly more general than existing data structures in the literature, encompassing, for example, a spectrum of problems classes from integer programming with indicator constraints to bilinear semidefinite programming. We also outline an automated rewriting system between equivalent formulations of a constraint. MathOptInterface has been implemented in practice, forming the foundation of a recent rewrite of JuMP, an open-source algebraic modeling language in the Julia language. The regularity of the MathOptInterface representation leads naturally to a general file format for mathematical optimization we call \textit{MathOptFormat}. In addition, the automated rewriting system provides modeling power to users while making it easy to connect new solvers to JuMP.}
}


\maketitle

\section{Introduction}

JuMP \citep{lubinComputingOperationsResearch2015,dunningJuMPModelingLanguage2017} is an algebraic modeling language for mathematical optimization written in the Julia language \citep{bezansonJuliaFreshApproach2017}. JuMP, like other algebraic modeling languages (e.g., AMPL \citep{fourer1990modeling}, Convex.jl \citep{convexjl}, CVX \citep{grant2014cvx}, CVXPY \citep{diamond2016cvxpy}, GAMS \citep{brook1988gams}, Pyomo \citep{hart2011pyomo,hart2017pyomo}, and YALMIP \citep{lofberg2004yalmip}), has an appearingly simple job: it takes a mathematical optimization problem written by a user, converts it into a \textit{standard form}, passes that standard form to a solver, waits for the solver to complete, then queries the solver for a solution and returns the solution to the user.

At the heart of this process is the definition of the standard form. By standard form, we mean a concrete data structure, specified either by an in-memory API or via a file format, of a mathematical optimization problem that the user and solver agree upon so that they can communicate. For example, based on the textbook presentation of linear programming (LP), one might assume that the following is a standard form accepted by solvers:
\begin{equation}\label{eq:lp_standard_form}
    \begin{array}{rl}
        \min\limits_{\mathbf{x} \in \mathbb{R}^N} & \mathbf{c}^\top \mathbf{x}\\
        \text{subject to:} & A\mathbf{x} = \mathbf{b}\\
                           & \mathbf{x} \ge \mathbf{0},
    \end{array}
\end{equation}
where $\mathbf{c}$ is a dense $N$-dimensional vector, $\mathbf{b}$ a dense $M$-dimensional vector, and $A$ a sparse\footnote{For this discussion, consider any standard sparse storage format, e.g., compressed sparse column or coordinate-list format.} $M\times N$ matrix.

Anyone who has interacted directly with LP solvers would know this is far from accurate. Some solvers allow linear constraints to have lower and upper bounds, so the user must pass $\mathbf{l} \le A\mathbf{x} \le \mathbf{u}$. Other solvers allow only one bound per row of $A$, but the user must also pass a vector of constraint senses (i.e., $=$, $\le$, or $\ge$), representing the problem constraints $A\mathbf{x} \triangleright \mathbf{b}$, where $\triangleright$ is the vector of constraint senses. Differences in these formulations also flow through to solutions, where, given an affine constraint $l \le \mathbf{a}^\top \mathbf{x} \le u$, some solvers may return two dual variables---one for each side of the constraint---whereas other solvers may return one dual variable corresponding to the active side of the constraint. In addition, some solvers may support variable bounds, whereas others based on conic methods may require variables to be non-negative.

Moreover, as mathematical optimization has matured, research has focused on formulating and solving new types of optimization problems. Each time a new type of objective function or constraint has been added to the modeler's toolbox, the standard form necessarily has had to change. For example, the commercial solvers MOSEK \citep{mosek} and Gurobi \citep{gurobi} have both developed independent (and incompatible) extensions to the MPS file format \citep{mps_ibm} to support quadratic objectives and constraints. \revisionb{This} has led to a fracturing of the optimization community as each sub-community developed a different standard form and solver for the problems of their interest. For example, nonlinear programming solvers often require the standard form:
\begin{equation*}
    \begin{array}{rl}
        \min\limits_{\mathbf{x} \in \mathbb{R}^N} & f(\mathbf{x})\\
        \text{subject to:} & g(\mathbf{x}) \le \mathbf{0}\\
                           & h(\mathbf{x}) = \mathbf{0},
    \end{array}
\end{equation*}
where $f: \mathbb{R}^N \mapsto \mathbb{R}$, $g: \mathbb{R}^N \mapsto \mathbb{R}^G$, and $h: \mathbb{R}^N \mapsto \mathbb{R}^H$ (and their respective derivatives) are specified via callbacks. In another sub-community, semidefinite solvers require the standard form:
\begin{equation*}
    \begin{array}{rl}
        \min\limits_{X \in \mathbb{R}^N\times\mathbb{R}^N} & \langle C, X \rangle\\
        \text{subject to:} & \langle A_i, X\rangle = b_i, \quad i=1,2,\ldots,M\\
                           & X \succeq 0,
    \end{array}
\end{equation*}
where $C$ and $A_i$ are $N\times N$ matrices, $b_i$ is a constant scalar, $\langle \cdot, \cdot\rangle$ denotes the inner product, and $X \succeq 0$ enforces the matrix $X$ to be positive semidefinite.

Even within communities, there can be equivalent formulations. For example, in conic optimization, solvers such as CSDP \citep{borchers1999csdp} accept what we term the \textit{standard conic} form:
\begin{equation}\label{eq:conic_standard_form}
    \begin{array}{rl}
        \min\limits_{\mathbf{x}\in\mathbb{R}^N} & \mathbf{c}^\top \mathbf{x}\\
        \text{subject to:} & A\mathbf{x} = \mathbf{b}\\
                           & \mathbf{x} \in \mathcal{K},
    \end{array}
\end{equation}
whereas others solvers, such as SCS \citep{scs}, accept what we term the \textit{geometric conic} form:
\begin{equation}\label{eq:conic_geometric_form}
    \begin{array}{rl}
        \min\limits_{\mathbf{x}\in\mathbb{R}^N} & \mathbf{c}^\top \mathbf{x}\\
        \text{subject to:} & A\mathbf{x} + \mathbf{b} \in \mathcal{K}\\
                           & \mathbf{x}\ \mathsf{free}.
    \end{array}
\end{equation}
Here, $\mathbf{c}$ is a $N$-dimensional vector, $A$ is an $M\times N$ matrix, $\mathbf{b}$ is an $M$-dimensional vector, and $\mathcal{K} \subseteq \mathbb{R}^N$ ($\mathcal{K}\subseteq\mathbb{R}^M$ for the geometric conic form) is a convex cone from some pre-defined list of supported cones.

\subsection{\revision{Contributions and outline}}

The variation in standard forms accepted by solvers makes writing a generalized algebraic modeling language such as JuMP difficult.
\revision{In this paper, we introduce three conceptual contributions to make this job easier:
\begin{enumerate}
    \item[(i)] We define \textit{MathOptInterface}, a new abstract data structure\footnote{By \textit{abstract} data structure, we mean that we omit discussions of details like which storage format to use for sparse matrices and how an API should be implemented in a specific programming language.} for representing mathematical optimization problems that generalizes the real-world diversity of the forms expected by solvers. 
    \item[(ii)] We describe an automated rewriting system \revisionb{based on} ``bridges'' between equivalent formulations of a constraint.
    \item[(iii)] We introduce a file format called MathOptFormat which is a direct serialization of MathOptInterface models into the JSON file format.
\end{enumerate}}

\revision{As a fourth contribution, we provide an implementation our ideas in the \texttt{MathOptInterface}\footnote{\url{https://github.com/jump-dev/MathOptInterface.jl}} library in Julia. This library is the foundation of a recent rewrite of JuMP. \texttt{MathOptInterface} was first released in February 2019, and provides a practical validation of our conceptual ideas. In addition, the implementation is a useful guide for others looking to implement our conceptual ideas in different programming languages.}

It is important to note that this paper deals with both the \textit{abstract} idea of the MathOptInterface standard form, and an implementation of this idea in Julia. To clearly distinguish between the two, we will always refer to the Julia constructs in \texttt{typewriter font}. Readers should note however, that our standard form is not restricted to the Julia language. Instead, it is intended to be a generic framework for thinking and reasoning about mathematical optimization. It is possible to write implementations in other languages such as Python; however, we chose Julia because JuMP is in Julia. Incidentally, the features of Julia make it well suited for implementing MathOptInterface in a performant way.

The rest of this paper is laid out as follows.
In Section \ref{sec:literature}, we review the approaches taken by other algebraic modeling languages and the old version of JuMP\revision{, before outlining the ideas behind the conception of MathOptInterface in Section \ref{sec:design_principles}.}
Then, in Section \ref{sec:standard_form}, we \revision{formally} introduce the MathOptInterface \revision{abstract data structure}, which is the main contribution of this paper.
In Section \ref{sec:bridges}, we \revision{introduce the automated constraint rewriting system}.
In Section \ref{sec:mathoptformat}, we present a new file format---called \textit{MathOptFormat}---for mathematical optimization that is based on MathOptInterface.
Finally, in Section \ref{sec:jump_1.0}, we describe how the introduction of \texttt{MathOptInterface} has influenced JuMP.

\section{Literature review}\label{sec:literature}

In this section, we review how existing modeling packages manage the conflict between the models provided by users and the standard forms expected by solvers. 

\subsection{A history of modeling packages}

\citet{orchard-haysHistoryMathematicalProgramming1984} provides a detailed early history of the inter-relationship between computing and optimization, beginning with the introduction of the simplex algorithm in 1947, through to the emergence of microcomputers in the early 1980s. Much of this early history was dominated by a punch-card input format called MPS \citep{mps_ibm}, which users created using problem-specific computer programs called \textit{matrix generators}.

However, as models kept getting larger, issues with matrix generators began to arise. \citet{fourer1983modeling} argues that the main issues were: (i) a lack of verifiability, which meant that bugs would creep into the matrix generating code; and (ii) a lack of documentation, which meant that it was often hard to discern what algebraic model the matrix generator actually produced. Instead of matrix generators, \citet{fourer1983modeling} advocated the adoption of \textit{algebraic modeling languages}. Algebraic modeling languages can be thought of as advanced matrix generators which build the model by parsing an algebraic model written in a human-readable domain-specific language. Two examples of early modeling languages that are still in wide-spread use are AMPL \citep{fourer1990modeling} and GAMS \citep{brook1988gams}.

Modeling languages allowed users to construct larger and more complicated models. In addition, they were extended to support the ability to model different types of programs, e.g., nonlinear programs, and specific constraints such as complementarity constraints. Because of the nonlinear constraints, modeling languages such as AMPL and GAMS represent constraints as \textit{expression graphs}. They are able to communicate with solvers through the \textit{NL} file-format \citep{gayWritingNlFiles2005}, or through modeling-language specific interfaces such as the \textit{AMPL Solver Library} \citep{gay1997hooking}. 

More recently, the last 15 years has seen the creation of modeling languages embedded in high-level programming languages. Examples include CVX \citep{grant2014cvx} and YALMIP \citep{lofberg2004yalmip} in MATLAB$^\text{\textregistered}$, CVXPY \citep{diamond2016cvxpy} and Pyomo \citep{hart2011pyomo} in Python, and Convex.jl \citep{convexjl} and JuMP in Julia.

Like AMPL and GAMS, Pyomo represents models using expression graphs and interfaces with solvers either through files (e.g., MPS and NL files), or, for a small number of solvers, via direct in-memory interfaces which convert the expression graph into each solver's specific standard form. 

YALMIP is a MATLAB$^\text{\textregistered}$-based modeling language for mixed-integer conic and nonlinear programming. It also has support for other problem classes, including geometric programs, parametric programs, and robust optimization. Internally, YALMIP represents conic programs in the geometric form~\eqref{eq:conic_geometric_form}.
Because of this design decision, if the user provides a model that is close to the standard conic form~\eqref{eq:conic_standard_form},
YALMIP must convert the problem back to the geometric form, introducing additional slack variables and constraints.
This can lead to sub-optimal formulations being passed to solvers.
As a work-around for this issue, YALMIP provides functionality for automatic dualizing conic models \citep{Lofberg2009}. Solving the dual instead of the primal can lead to significantly better performance in some cases; however, the choice of when to dualize the model is left to the user.

Since YALMIP represents the problem in geometric form and has no interface allowing users to specify the cones the variables belong to, when dualizing it needs to reinterpret the affine constraints when the affine expression contains only one variable as a specification of the cone for the variable.
Given that there is no unique interpretation of the dual form if a variable is interpreted to belong to two or more cones, in such cases YALMIP considers one of the constraints as a variable constraint and the others as affine constraints.
As described by \citet{Lofberg2009}, YALMIP uses the following heuristic if the constraints are on different cones:
``Notice that it is important to detect the primal cone variables in a certain order, starting with SDP cones, then SOCP cones, and finally LP cones.''

The rigidity of a standard form chosen by a modeling language such as YALMIP also limits the structure that the user is able to transmit to the solver.
The resulting transformations needed to make a problem fit in the standard form can have significant impacts on the runtime performance of the solver.
For example, in addition to formulation~\eqref{eq:conic_standard_form}, SDPNAL+ \citep{sun2019sdpnal+} supports adding bounds on the variables and adding affine interval constraints.
Forcing the problem to fit in the standard form~\eqref{eq:conic_standard_form} requires the addition of slack variables and equality constraints that have a negative impact on the performance of SDPNAL+
as described in \citet{sun2019sdpnal+}:
\begin{quote}
    ``The final number of equality constraints present in the data input to SDPNAL+ can also be substantially fewer than those present in \revisionb{[...]} YALMIP. It is important to note here that the number of equality constraints present in the generated problem data can greatly affect the computational efficiency of the solvers, especially for interior-point based solvers.''
\end{quote}

CVXPY \citep{diamond2016cvxpy} is a modeling language for convex optimization in Python. A notable feature of CVXPY is that it is based on disciplined convex programming \citep{grant2006disciplined}; this is a key difference from many other modeling languages, including JuMP. The rules of disciplined convex programming mean that the convexity of a user-provided function can be inferred at construction time in an axiomatic way. This has numerous benefits for computation, but restricts the user to formulating models with a reduced set of operators (called \textit{atoms}), for which the convexity of the atom is known. \revisionb{Recent upgrades to CVXPY have added support for log-log convex programming and quasi-convex programming \citep{agrawal2019disciplined}. CVXPY also supports mixed-integer programs.}

One feature of CVXPY that the re-write of JuMP does inherit is the concept of a \textit{reduction}. A reduction is a transformation of one problem into an equivalent form. Reductions allow CVXPY to re-write models formulated by the user into equivalent models that solvers accept \citep{cvxpy_rewriting}. Examples of reductions implemented in CVXPY include \texttt{Complex2Real}, which lifts complex-valued variables into the real domain by introducing variables for the real and imaginary terms, and \revisionb{\texttt{Dgp2Dcp}, which converts a disciplined geometric program into a disciplined convex program \citep{agrawal2019disciplined}}. Reductions can be chained together to form a sequence of reductions. CVXPY uses pre-defined chains of reductions to convert problems from the form given by the user into a standard form required by a solver. 

\revision{Constraint transformations have also been explored in the context of the constraint programming language MiniZinc \citep{nethercote2007minizinc}. MiniZinc is a standalone modeling language similar to AMPL and GAMS. To communicate with solvers, MiniZinc compiles problems formulated by the user into a low-level file format called FlatZinc \citep{marriott2008design}. During this compilation step, the user's model is rewritten into a form supported by the targeted solver. In particular, MiniZinc  allows users to write constraints such as \texttt{alldifferent}$(\mathbf{x})$, which enforces that no two elements in the vector $\mathbf{x}$ can take the same value. These constraints can be either passed directly to constraint programming solvers, or reformulated into a mixed-integer linear program by including a redefinition file in the model's source code. MiniZinc provides a default library of redefinition files which can be chosen by the user. However, if the user has advanced knowledge of their problem and solver, they can write a new definition that will be used in-place of the default transform \citep{brand2008flexible,belov2016improved}.}

\revision{Like CVXPY's reductions, MiniZinc's transformations need to be chosen ahead-of-time. A key innovation in the bridging system we describe in Section \ref{sec:bridges} is that the chain of transformations are automatically chosen as the model is built at \textit{run-time}.  However, JuMP's reductions apply only to classes of constraints, variables, and objective functions, rather than applying global transformations as CVXPY does.}

\subsection{A history of file formats}\label{sec:file_format_history}

In Section \ref{sec:mathoptformat}, we introduce a new file format for mathematical optimization. \revision{Given the existing, widely adopted file formats that have served the community well to date, creating a new format is not a decision that we made lightly.} Our main motivation was to create a way to serialize MathOptInterface problems to disk. However, to fully understand our justification for creating a new format, it is necessary to give a brief history of the evolution of file formats in mathematical optimization.

As we have outlined, in order to use an optimization solver, it is necessary to communicate a model instance to the solver. This can be done either through an in-memory interface, or through a file written to disk. File formats are also used to collate models into instance libraries for benchmarking purposes, e.g., CBLIB \citep{fribergCBLIB2014Benchmark2016}, \revision{MINLPLib \citep{bussieck2003minlplib,minlplib}}, and MIPLIB \citep{miplib2017}.

Many different instance formats have been proposed over the years, but only a few (such as MPS \citep{mps_ibm}) have become the industry standard. Each format is a product of its time in history and the problem class it tried to address. For example, we retain the rigid input format of the MPS file that was designed for 1960s punch-cards despite the obsolescence of this technology \citep{orchard-haysHistoryMathematicalProgramming1984}. Although the MPS format has since been extended to problem classes such as nonlinear and stochastic linear programming, MPS was not designed with extensibility in mind. This has led some authors (e.g., \citet{fribergCBLIB2014Benchmark2016}) to conclude that developing a new format is easier than extending the existing MPS format.

The LP file-format is an alternative to the MPS file-format that is human-readable and row-oriented \citep{lp_format}. However, there is no longer a single standard for the LP file-format. This has led to subtle differences between implementations in different readers that hampers the usefulness of the format as a medium for interchange. Much like the MPS file, the LP file is also limited in the types of problems it can represent and was not designed for extensibility.

In contrast to the LP file, the~NL file \citep{gayWritingNlFiles2005} explicitly aims for machine-readability at the expense of human-readability. It is also considerably more flexible in the problem classes it can represent (in particular, nonlinear functions are supported). However, once again, the format is not extensible to new problem formats and has limited support for conic problems.

\revision{GAMS scalar format \citep{bussieck2003minlplib}, is a GAMS-based file format for serializing nonlinear programs. The GAMS scalar format uses a subset of the full GAMS syntax, and so it is human readable. However, it has limited support for conic programs. That is, simple second-order cone constraints can be specified, although this feature has been deprecated.}

The OSiL format \citep{fourerOSiLInstanceLanguage2010} is an XML-based file format that targets a broad range of problem classes. In developing OSiL, Fourer et al.~identified many of the challenges and limitations of previous formats and attempted to overcome them. In particular, they choose to use XML as the basis for their format to remove the burden of writing custom readers and writers for each programming language that wished to interface with optimization software, allowing more focus on the underlying data structures. XML is also human-readable and can be rigidly specified with a schema to prevent the proliferation of similar, but incompatible versions. The XML approach has been extended to support multiple problem classes including nonlinear, stochastic, and conic.

However, despite the many apparent advantages of the OSiL format, we believe it has enough short-comings to justify the development of a new instance format. The main reason is the lack of a strong, extensible standard form. A secondary reason is the waning popularity of XML in favor of simpler formats such as JSON.

\subsection{The previous design of JuMP}

Until recently, JuMP, described by \citet{lubinComputingOperationsResearch2015} and \citet{dunningJuMPModelingLanguage2017}, featured a collection of three different standard forms: (i) a \textit{linear-quadratic} standard form for specifying problems with linear, quadratic, and integrality constraints; (ii) a \textit{conic} standard form for specifying problems with linear, conic, and integrality constraints; and (iii) a \textit{nonlinear} standard form for specifying problems with nonlinear constraints.

In code, the three standard forms were implemented in an intermediate layer called \texttt{MathProgBase}. As the first step, JuMP converted the problem given by the user into one of the three \texttt{MathProgBase} standard forms. Underneath \texttt{MathProgBase}, each solver required a thin layer of Julia code (called \textit{wrappers}) that connected the solver's native interface (typically written in C) to one of the three standard forms, e.g, Clp \citep{clp} to \textit{linear-quadratic}, SCS \revision{\citep{scs}} to \textit{conic}, and Ipopt \citep{wachter2006implementation} to \textit{nonlinear}. 
There were also automatic converters between the standard forms, including linear-quadratic to conic, conic to linear-quadratic, and linear-quadratic to nonlinear. This enabled, for example, the user to formulate a linear program and solve it with Ipopt (a nonlinear programming solver). Figure \ref{fig:jump_architecture} visualizes the architecture of JuMP just described.

\begin{figure}[!ht]
    \centering
    \resizebox{0.6\textwidth}{!}{%
    \begin{tikzpicture}[align=center]
        \node[my_node] (jump)       at (4, 6) {JuMP};
        \node[mpb_node, dashed] (mpb) at (4, 3) {};
        \node[my_node] (lq)         at (4, 3) {Linear-Quadratic\\standard form};
        \node[my_node] (con)        at (0, 3) {Conic\\standard form};
        \node[my_node] (nl)         at (8, 3) {Nonlinear\\standard form};
        \node[my_node] (lq_solver)  at (4, 0) {Linear-Quadratic\\solver};
        \node[my_node] (con_solver) at (0, 0) {Conic\\solver};
        \node[my_node] (nl_solver)  at (8, 0) {Nonlinear\\solver};
        \node[] (mpb_text)  at (-0.5, 4.65) {\texttt{MathProgBase}};
    
        \path[very thick]
        	(jump)	edge [] node [] {} (lq)
        	        edge [] node [] {} (con)
        	        edge [] node [] {} (nl)
        	(lq)	edge [] node [] {} (lq_solver)
        	        edge [] node [] {} (con)
        	        edge [] node [] {} (nl)
        	        edge [] node [] {} (jump)
        	(con)	edge [] node [] {} (con_solver)
        	        edge [] node [] {} (lq)
        	        edge [] node [] {} (jump)
        	(nl)	edge [] node [] {} (nl_solver)
        	        edge [] node [] {} (jump)
        	(lq_solver)	edge [] node [] {} (lq)
        	(con_solver)	edge [] node [] {} (con)
        	(nl_solver)	edge [] node [] {} (nl);
    \end{tikzpicture}
    }
    \caption{Architecture of JuMP before it switched to MathOptInterface. JuMP models are communicated through to solvers via the \texttt{MathProgBase} interface layer (dashed line), which consists of three standard forms.}
    \label{fig:jump_architecture}
\end{figure}

The design of \texttt{MathProgBase} made it impossible to flexibly combine different standard forms. For example, JuMP could not communicate a program with both nonlinear and conic constraints to a solver.

\section{\revisionb{Design principles}}\label{sec:design_principles}

\revisionb{In mid-2017, we decided to move away from \texttt{MathProgBase} and re-write JuMP based on a new, all-encompassing standard form and interface layer that became the MathOptInterface specification and the \texttt{MathOptInterface} package. The main reasons for this move were the inflexibility of \texttt{MathProgBase} for mixing problem classes, the difficulty of adding support for new types of constraints, and various other early design decisions that became hard to change\footnote{Although out of scope for this paper, these shortcomings were addressed by the features briefly mentioned in Section~\ref{sec:other_features}.}.}


\revision{On one hand, the intended scope of the new abstraction was relatively narrow. At a minimum, MathOptInterface needed to support all the problem classes that JuMP supported at the time (via \texttt{MathProgBase}) and combinations thereof. We did not try to cover new paradigms like multi-level or stochastic optimization. Additionally, because of the complexities surrounding automatic differentiation, we deferred first-class support for nonlinear programming for future work. Finally, we were willing to accept some trade-off in performance for simplicity of the design, although we have taken care to ensure that the abstraction layer does not induce a bottleneck in the solve process.}

\revisionb{On the other hand, the experience writing of \texttt{MathProgBase} led us to two important design principles that underlie the aspects of MathOptInterface covered in this paper:
\begin{itemize}
    \item[(i)] MathOptInterface should be \textit{extensible} to new types of constraints.
    \item[(ii)] MathOptInterface should expose the tension that exists between the generality of how users expect to express models, and the rigidity of how solvers expect to receive models.
\end{itemize}}

\revisionb{With regard to extensibility,} \revision{we had learned from JuMP development under the former \texttt{MathProgBase} interface that the set of possible constraints that users want to model can overwhelm our ability as a small team of open-source developers to accommodate. For example, we did not manage to support indicator or complementarity constraints in \texttt{MathProgBase} because doing so would have required simultaneous invasive changes in \texttt{MathProgBase} and JuMP, a task effectively too large for even a committed contributor who was not part of the core team. A goal for \texttt{MathOptInterface} was to have a well-documented and accessible structure for introducing new constraint types that required few, if any, changes to JuMP. We settled on a very \textit{regular} representation for constraints (and for other aspects of the abstraction), so much so that JuMP can process new types of constraints defined in add-on packages.}

\revisionb{The second} \revision{important consideration was the idea to expose the tension we had observed between how users expect to express their models' constraints and how solvers accept the constraints. This could be, for example, the difference between a typical LP model and the more rigid standard conic form~\eqref{eq:conic_standard_form} that doesn't support explicit bounds on variables, the difference between the second order cone and the rotated second order cone (each is a linear transformation of the other), or the difference between indicator constraints and pure integer programming formulations (e.g., big-M). We wanted modelers and solvers to speak the same ``language,'' so that solvers can advertise which constraints they support, modelers have an array of options for how to express their model, and the two sides can be bridged by transformations of a common data structure. It was intended that this bridging could happen, either (i) by an automated rewriting system, as we later describe, or (ii) by modelers deciding to rewrite their model in a format closer to what a solver natively supports for additional control and performance, or (iii) by solver authors deciding to support new types of constraints in response to demand from users. The latter two goals reflect a view of modeling and solver development as dynamic processes with multiple self-motivated agents. Indeed, while this tension between modelers and solvers is always visible to developers of algebraic modeling interfaces, our idea was to expose it more concretely and programmatically so that any motivated modeler or solver developer could take their own steps to address it while remaining within a common abstraction instead of reverting to solver-specific interfaces.}

\section{MathOptInterface}\label{sec:standard_form}

\revisionb{The considerations in the previous section guided our design of} MathOptInterface, an abstract specification for a data structure for mathematical optimization problems\revision{, which we now describe formally. MathOptInterface} represents problems in the following form:
\begin{equation}\label{eq:moi_standard_form}
    \begin{array}{rl}
        \min\limits_{\mathbf{x} \in \mathbb{R}^N} & f_0(\mathbf{x})\\
        \text{subject to:} & f_i(\mathbf{x}) \in S_i,\quad i=1,2,\ldots,I,
    \end{array}
\end{equation}
with \textit{functions} $f_i:\;\mathbb{R}^N \mapsto \mathbb{R}^{M_i}$ and \textit{sets} $S_i \subseteq \mathbb{R}^{M_i}$ drawn from a pre-defined set of functions $\mathcal{F}$ and sets $\mathcal{S}$. 

The sets $\mathcal{F}$ and $\mathcal{S}$ are provided in Section \ref{sec:moi_functions} (for $\mathcal{F}$) and the Online Supplement (for $\mathcal{S}$). In addition, we provide a concrete description in the form of a JSON schema \citep{json_schema} as part of the \textit{MathOptFormat} file format described in Section \ref{sec:mathoptformat}. In the JuMP ecosystem, the definitions of supported functions and sets are contained in the Julia package \texttt{MathOptInterface}. JuMP additionally allows third-party packages to extend the set of recognized functions and sets $\mathcal{F}$ and $\mathcal{S}$ at run-time, but this is not required for an implementation of MathOptInterface; for example, MathOptFormat does not allow extensions.

Since constraints are formed by the combination of a function and a set, we will often refer to constraints by their \textit{function}-in-\textit{set} pairs. The key insight is the ability to mix-and-match a small number of pre-defined functions and sets to create a wide variety of different problem classes.

We believe model~\eqref{eq:moi_standard_form} is very general, and encompasses almost all of existing deterministic mathematical optimization with real-valued variables. (An extension to complex numbers could be achieved by replacing $\mathbb{R}$ with $\mathbb{C}$.) Readers should note that when the objective is vector-valued, the objective vectors are implicitly ranked according to  partial ordering such that if $\mathbf{y_1} = f_0(\mathbf{x_1})$ and $\mathbf{y_2} = f_0(\mathbf{x_2})$, then $\mathbf{y_1} \le \mathbf{y_2} \iff \mathbf{y_2} - \mathbf{y_1} \in \mathbb{R}_+^{M_0}$. In the future, we plan to extend MathOptInterface to model general vector-valued programs, which define the partial ordering in terms of a convex cone $\mathcal{C}$ (see, e.g., \citet{lohne2011vector}). However, we omit a full description of this extension because we do not have a practical implementation.

\subsection{Functions}\label{sec:moi_functions}

MathOptInterface defines the following functions in the set $\mathcal{F}$:
\begin{itemize}
    \item The \textsf{SingleVariable} function \revision{$f:\mathbb{R}^N\mapsto \mathbb{R}$ with} $f(\mathbf{x}) = \mathbf{e_i}^\top \mathbf{x}$, where $\mathbf{e_i}$ is an $N$-dimensional vector of zeros with a $1$ in the $i^{th}$ element.
    \item The \textsf{VectorOfVariables} function \revision{$f:\mathbb{R}^N\mapsto \mathbb{R}^M$ with} $f(\mathbf{x}) = [x_{i_1}, x_{i_2}, \ldots, x_{i_M}]$, where $i_j \in \{1,2,\ldots, N\}$ for all $j \in 1,\ldots,M$.
    \item The \textsf{ScalarAffineFunction} \revision{$f:\mathbb{R}^N\mapsto \mathbb{R}$ with} $f(\mathbf{x}) = \mathbf{a}^\top \mathbf{x} + b$, where $a$ is a sparse $N$-dimensional vector and $b$ is a scalar constant.
    \item The \textsf{VectorAffineFunction} \revision{$f:\mathbb{R}^N\mapsto \mathbb{R}^M$ with} $f(\mathbf{x}) = A^\top \mathbf{x} + \mathbf{b}$, where $A$ is a sparse $M\times N$ matrix and $\mathbf{b}$ is a dense $M$-dimensional vector.
    \item The \textsf{ScalarQuadraticFunction} \revision{$f:\mathbb{R}^N\mapsto \mathbb{R}$ with} $f(\mathbf{x}) = \frac{1}{2}\mathbf{x}^\top Q \mathbf{x} + \mathbf{a}^\top \mathbf{x} + b$, where $Q$ is a sparse $N\times N$ matrix, $\mathbf{a}$ is a sparse $N$-dimensional vector, and $b$ is a scalar constant.
    \item The \textsf{VectorQuadraticFunction} \revision{$f:\mathbb{R}^N\mapsto \mathbb{R}^M$ with} $f(\mathbf{x}) = \left[\mathbf{x}^\top Q_1 \mathbf{x},\cdots, \mathbf{x}^\top Q_i \mathbf{x},\cdots, \mathbf{x}^\top Q_M \mathbf{x}\right]^\top + A^\top \mathbf{x} + \mathbf{b}$, where $Q_i$ is a sparse $N\times N$ matrix for $i=1,2,\ldots M$, $A$ is a sparse $M\times N$ matrix, and $\mathbf{b}$ is a dense $M$-dimensional vector.
\end{itemize}

Notably missing from this list is \textsf{ScalarNonlinearFunction} and \textsf{VectorNonlinearFunction}. At present, \texttt{MathOptInterface} defines (\revision{from the legacy of \texttt{MathProgBase}}) a separate mechanism for declaring a block of nonlinear constraints $\mathbf{l} \le g(\mathbf{x}) \le \mathbf{u}$ and/or a nonlinear objective $f(\mathbf{x})$. \revision{Integrating nonlinear functions into \texttt{MathOptInterface} as a first-class object will likely require further generalizations, particularly around modular support for automatic differentiation. We leave such details to future work.}

Moreover, note that many of the function definitions are redundant, e.g., a \textsf{ScalarAffineFunction} is a \textsf{VectorAffineFunction} where $M=1$. The reason for this redundancy \revision{is to expose the variety of ways modelers and solvers prefer to express their problems, as discussed in Section~\ref{sec:design_principles}}.

\subsection{Sets}\label{sec:moi_sets}

The set of sets supported by MathOptInterface, $\mathcal{S}$, contains a large number of elements. The complete list is given in the Online Supplement and by the JSON schema described in Section \ref{sec:mathoptformat}, so we only present some of the more common sets that will later be referenced in this paper:

\begin{itemize}
    \item The \textsf{LessThan} set $(-\infty, u]$ where $u\in\mathbb{R}$
    \item The \textsf{GreaterThan} set $[l, \infty)$ where $l\in\mathbb{R}$
    \item The \textsf{Interval} set $[l, u]$, where $l\in\mathbb{R}$ and $u\in\mathbb{R}$
    \item The \textsf{Integer} set $\mathbb{Z}$
    \item The \textsf{Nonnegatives} set $\left\{\mathbf{x}\in\mathbb{R}^N\;:\;\mathbf{x} \ge \mathbf{0}\right\}$
    \item The \textsf{Zeros} set $\{\mathbf{0}\} \subset \mathbb{R}^N$
    \item The \textsf{SecondOrderCone} set $\left\{(t, \mathbf{x})\in\mathbb{R}^{1+N}\;:\; ||\mathbf{x}||_2 \le t\right\}$\revision{, where $\mathbf{x}\in\mathbb{R}^N$}
    \item The \textsf{RotatedSecondOrderCone} set $\left\{(t, u, \mathbf{x})\in\mathbb{R}^{2+N}\;:\; ||\mathbf{x}||_2^2 \le 2tu, t\ge 0, u \ge 0\right\}$\revision{, where $\mathbf{x}\in\mathbb{R}^N$}
\end{itemize}

MathOptInterface also defines sets like the positive semidefinite cone, and even sets that are not cones or standard sets like \textsf{Interval} and \textsf{Integer}. For example, MathOptInterface defines the \textsf{SOS1} and \textsf{SOS2} sets, which are special ordered sets of Type I and Type II respectively \citep{beale1970special}. In addition, it also defines the \textsf{Complements} set, which can be used to specify mixed complementarity constraints \citep{dirkse1995path}. See the Online Supplement for more details.

To demonstrate how these functions and sets can be combined to create mathematical programs, we now consider a number of examples.

\subsection{Example: linear programming}

Linear programs are often given in the form:
\begin{subequations}\begin{align}
    \min\limits_{\mathbf{x}\in\mathbb{R}^N}\ & \mathbf{c}^\top \mathbf{x} + 0 \label{eq:lp1_obj}\\
    \text{subject to:}\ & A\mathbf{x} \ge \mathbf{b}, \label{eq:lp1_con}
\end{align}\end{subequations}
where $\mathbf{c}$ is a vector with $N$ elements, $A$ is an $M\times N$ matrix, and $\mathbf{b}$ is a vector with $M$ elements. 

In the MathOptInterface standard form, \revision{objective~\eqref{eq:lp1_obj}} is the \textsf{ScalarAffineFunction} $f_0(\mathbf{x}) = \mathbf{c}^\top \mathbf{x} + 0$, \revision{and constraint~\eqref{eq:lp1_con} is composed of} the \textsf{VectorAffineFunction} $f_1(\mathbf{x}) = A\mathbf{x} - \mathbf{b}$ and the $M$-dimensional \textsf{Nonnegatives} set.

\subsection{Example: multi-objective problems with conic constraints}

Because of its generality, MathOptInterface is able to represent problems that do not neatly fit into typical standard forms. For example, here is a multi-objective mathematical program with a second-order cone constraint:
\begin{subequations}\begin{align}
    \min\limits_{\mathbf{x} \in \mathbb{R}^N}\ & C \mathbf{x} + \mathbf{0} \label{eq:mo_obj}\\
    \text{subject to:}\ & A\mathbf{x} = \mathbf{b} \label{eq:mo_con1}\\
                       & ||x_2, x_3, \ldots, x_N||_2 \le x_1 \label{eq:mo_con2}\\
                       & l_i \le x_i \le u_i, i=1,2,\ldots, N, \label{eq:mo_con3}
\end{align}\end{subequations}
where $C$ is an $P\times N$ matrix, $A$ is an $M \times N$ matrix, $\mathbf{b}$ is an $M$-dimensional vector, and $l_i$ and $u_i$ are constant scalars.

\revision{In the MathOptInterface standard form, objective~\eqref{eq:mo_obj} is the \textsf{VectorAffineFunction} $f_0(\mathbf{x}) = C\mathbf{x} + \mathbf{0}$, constraint~\eqref{eq:mo_con1} is composed of the \textsf{VectorAffineFunction} $f_1(\mathbf{x}) = A\mathbf{x} - \mathbf{b}$ and the $M$-dimensional \textsf{Zeros} set, constraint~\eqref{eq:mo_con2} is composed of the \textsf{VectorOfVariables} function $f_2(\mathbf{x}) = \mathbf{x}$ and the \textsf{SecondOrderCone} set, and constraints~\eqref{eq:mo_con3} are composed of the \textsf{SingleVariable} functions $f_{2 + i}(\mathbf{x}) = x_i$ and the \textsf{Interval} sets $[l_i, u_i]$.}

\subsection{Example: special ordered sets}

Many mixed-integer solvers support a constraint called a \textit{special ordered set} \citep{beale1970special}. There are two types of special ordered sets. Special ordered sets of type I require that at most one variable in an ordered set of variables can be non-zero. Special ordered sets of type II require that at most two variables in an ordered set of variables can be non-zero, and, if two variables are non-zero, they must be adjacent in the ordering.

An example of a problem with a special ordered set of type II constraint is as follows:
\begin{subequations}\begin{align}
    \min\limits_{\mathbf{x} \in \mathbb{R}^N}\ & \mathbf{c}^\top \mathbf{x} + 0 \label{eq:sos_obj}\\
    \text{subject to:}\ & A\mathbf{x} \ge \mathbf{b} \label{eq:sos_con1}\\
                       & [x_1, x_2, x_3] \in \textsf{SOS}_{\textsc{II}}([1, 3, 2]). \label{eq:sos_con2}
\end{align}\end{subequations}
Here, the weights on the variables imply an ordering $x_1$, $x_3$, $x_2$.

\revision{In the MathOptInterface standard form, objective~\eqref{eq:sos_obj} is the \textsf{ScalarAffineFunction} $f_0(\mathbf{x}) = \mathbf{c}^\top \mathbf{x} + 0$, constraint~\eqref{eq:sos_con1} is composed of the \textsf{VectorAffineFunction} $f_1(\mathbf{x}) = A\mathbf{x} - \mathbf{b}$ and the $M$-dimensional \textsf{Nonnegatives} set, and constraint~\eqref{eq:sos_con2} is composed of the \textsf{VectorOfVariables} function $f_{2}(\mathbf{x}) = [x_1, x_2, x_3]$ and the \textsf{SOS2} set \textsf{SOS2([1, 3, 2])}.}

\subsection{Example: mixed-complementarity problems}

A mixed-complementarity problem, which can be solved by solvers such as PATH \citep{dirkse1995path}, can be defined as follows:
\begin{equation}\label{eq:mcp}
    \begin{array}{rl}
        \min\limits_{\mathbf{x} \in \mathbb{R}^N}\ & 0\\
        \text{subject to:} & A\mathbf{x} + \mathbf{b} \perp \mathbf{x}\\
                           & \mathbf{l} \le \mathbf{x} \le \mathbf{u},
    \end{array}
\end{equation}
where $A$ is an $N\times N$ matrix, and $\mathbf{b}$, $\mathbf{l}$, and $\mathbf{u}$ are $N$-dimensional vectors. The constraint $A\mathbf{x} + \mathbf{b} \perp \mathbf{x}$ requires that the following conditions hold in an optimal solution:
\begin{itemize}
    \item if $x_i = l_i$, then $\mathbf{e_i}^\top A\mathbf{x} + b_i \ge 0$;
    \item if $l_i < x_i < u_i$, then $\mathbf{e_i}^\top A\mathbf{x} + b_i = 0$; and
    \item if $x_i = u_i$, then $\mathbf{e_i}^\top A\mathbf{x} + b_i \le 0$.
\end{itemize}

Thus, we can represent model~\eqref{eq:mcp} in the MathOptInterface standard form as:
\begin{subequations}\begin{align}
\min\limits_{\mathbf{x} \in \mathbb{R}^N} & 0 \label{eq:mcp_obj}\\
        \text{subject to:} & \begin{bmatrix}A \\ I\end{bmatrix}\mathbf{x} + \begin{bmatrix}\mathbf{b} \\ \mathbf{0}\end{bmatrix} \in \textsf{Complements}() \label{eq:mcp_con1}\\
                           & x_i \in [l_i, u_i], i=1,2,\ldots,N. \label{eq:mcp_con2}
\end{align}\end{subequations}

\revision{Here, objective~\eqref{eq:mcp_obj} is the \textsf{ScalarAffineFunction} $f_0(\mathbf{x}) = 0$, constraint~\eqref{eq:mcp_con1} is composed of the \textsf{VectorAffineFunction} $f_1(x) = [A; I]\mathbf{x} + [\mathbf{b}; \mathbf{0}]$ and the $N$-dimensional \textsf{Complements} set, and constraints~\eqref{eq:mcp_con2} are composed of the \textsf{SingleVariable} functions $f_{1 + i}(\mathbf{x}) = x_i$ and \textsf{Interval} sets $[l_i, u_i]$.}

\section{Bridges}\label{sec:bridges}

\revision{By defining a small list of functions and sets, we obtain a large number of different constraint types. This design is naturally extensible and captures the diverse ways that solvers natively accept input, satisfying our two main design considerations. However, by creating many ways to express constraints, we also have to confront the challenge of translating between mathematically equivalent forms.} 

For example, the constraint $l \le \mathbf{a}^\top \mathbf{x} \le u$ can be formulated in many ways, three of which are listed here:
\begin{itemize}
    \item using the original formulation: $l \le \mathbf{a}^\top \mathbf{x} \le u$ (\textsf{ScalarAffineFunction}-in-\textsf{Interval});
    \item by splitting the constraint into $\mathbf{a}^\top \mathbf{x} \le u$ (\textsf{ScalarAffineFunction}-in-\textsf{LessThan}) and $\mathbf{a}^\top \mathbf{x} \ge l$ (\textsf{ScalarAffineFunction}-in-\textsf{GreaterThan}); or
    \item by introducing a slack variable $y$, with constraints $\mathbf{a}^\top \mathbf{x} - y = 0$ (\textsf{ScalarAffineFunction}-in-\textsf{EqualTo}) and $l \le y \le u$ (\textsf{SingleVariable}-in-\textsf{Interval}).
\end{itemize}
This generality means that solver authors need to decide which functions and sets to support, users need to decide how to formulate constraints, and modeling language developers need to decide how to translate constraints between the user and the solver.

One approach to this problem is to require every solver to implement an interface to every combination of function-in-set that the user could provide, and inside each solver transform the user-provided constraint into a form that the solver natively understands. However, as the number of functions and sets increases, this approach quickly becomes burdensome.

An alternative approach, and the one implemented in \texttt{MathOptInterface}, is to centralize the problem transformations into a collection of what we call \textit{bridges}. A bridge is a thin layer that transforms a function-in-set pair into an equivalent list of function-in-set pairs. An example is the transformation of a \textsf{ScalarAffineFunction}-in-\textsf{Interval} constraint into a \textsf{ScalarAffineFunction}-in-\textsf{LessThan} and a \textsf{ScalarAffineFunction}-in-\textsf{GreaterThan} constraint. The bridge is also responsible for reversing the transform to provide information such as dual variables back to the user.

\texttt{MathOptInterface} defines a large number of bridges. For example, there are \texttt{slack} bridges, which convert inequality constraints like $\mathbf{a}^\top \mathbf{x} \ge b$ into equality constraints like $\mathbf{a}^\top \mathbf{x} - y = b$ by adding a slack variable $y \ge 0$.
There are also bridges to convert between different cones. For example there is a bridge to convert a rotated second-order cone into a second-order cone using the following relationship:
$$2tu \ge || \mathbf{x} ||_2^2 \iff (t/\sqrt 2 + u/\sqrt 2)^2 \ge || \mathbf{x} ||_2^2 + (t/\sqrt 2 - u/\sqrt2)^2.$$

Bridges can also be nested to allow multiple transformations. For example, a solver that supports only \textsf{ScalarQuadraticFunction}-in-\textsf{EqualTo} constraints can support \textsf{RotatedSecondOrderCone} constraints via a transformation into a \textsf{SecondOrderCone} constraint, then into a \textsf{ScalarQuadraticFunction}-in-\textsf{LessThan} constraint, and then into a \textsf{ScalarQuadraticFunction}-in-\textsf{EqualTo} via a \texttt{slack} bridge.

The proliferation in the number of these bridges leads to a new challenge: there are now multiple ways of transforming one constraint into an equivalent set of constraints via chains of bridges. To demonstrate this, consider bridging a \textsf{ScalarAffineFunction}-in-\textsf{Interval} constraint into a form supported by a solver which supports only \textsf{SingleVariable}-in-\textsf{GreaterThan} and \textsf{ScalarAffineFunction}-in-\textsf{EqualTo} constraints. Two possible reformulations are given in Figure \ref{fig:bridging_scalar_affine_example}. 

\begin{figure}[!ht]
\centering
\resizebox{0.9\textwidth}{!}{
\begin{tikzpicture}[align=center]
    \node[draw, ellipse] (a) at (-2,0) {$l \le \mathbf{a}^\top \mathbf{x} \le u$};
    \node[] (b) at  (3,1) {$l \le \mathbf{a}^\top \mathbf{x}$};
    \node[] (c) at  (3,-1) {$\mathbf{a}^\top \mathbf{x} \le u$};
    \node[draw] (d) at (8,1.5) {$\mathbf{a}^\top \mathbf{x} - y = l$};
    \node[draw] (e) at (7.4,0.5) {$y \ge 0$};
    \node[draw] (f) at (8,-0.5) {$\mathbf{a}^\top \mathbf{x} + y = u$};
    \node[draw] (g) at (7.4,-1.5) {$y \ge 0$};
    \node[] (b1) at  (2.3,0) {\revision{\texttt{split-interval}}};
    \node[] (b2) at  (6,1) {\revision{\texttt{slack}}};
    \node[] (b3) at  (6,-1) {\revision{\texttt{slack}}};
    \path[->, out=0, in=180]
    	(a) edge (b)
    	    edge (c)
    	(b) edge (d)
    	    edge (e)
    	(c) edge (f)
    	    edge (g);
    \node[draw, ellipse] (a2) at  (-2,-4.5) {$l \le \mathbf{a}^\top \mathbf{x} \le u$};
    \node[] (b2) at  (3,-3.5) {$l \le \mathbf{a}^\top \mathbf{x}$};
    \node[] (c2) at  (3,-5.5) {$\mathbf{a}^\top \mathbf{x} \le u$};
    \node[draw] (d2) at  (8,-3) {$\mathbf{a}^\top \mathbf{x} - y = l$};
    \node[draw] (e2) at  (7.4,-4) {$y \ge 0$};
    \node[] (f2) at  (8,-5.5) {$-u \le -\mathbf{a}^\top \mathbf{x}$};
    \node[draw] (g2) at  (13.4,-5) {$-\mathbf{a}^\top \mathbf{x} - y = -u$};
    \node[draw] (h2) at  (12.4,-6) {$y \ge 0$};
    \node[] (b4) at  (2.3,-4.5) {\revision{\texttt{split-interval}}};
    \node[] (b5) at  (6, -3.5) {\revision{\texttt{slack}}};
    \node[] (b6) at  (11,-5.5) {\revision{\texttt{slack}}};
    \node[] (b7) at  (5.5,-6) {\revision{\texttt{flip-sign}}};
    \path[->, out=0, in=180]
    	(a2) edge (b2)
    	     edge (c2)
    	(b2) edge (d2)
    	     edge (e2)
    	(c2) edge (f2)
    	(f2) edge (g2)
    	     edge (h2);
\end{tikzpicture}}
\caption{Two equivalent solutions to the problem of bridging a \textsf{ScalarAffineFunction}-in-\textsf{Interval} constraint (oval nodes). Outlined rectangular nodes represent constraint actually added to the model. Nodes with no outline are intermediate nodes. \revision{\texttt{Typewriter} font describes the bridge used in the transformation.}}
\label{fig:bridging_scalar_affine_example}
\end{figure}

The first reformulation converts $l \le \mathbf{a}^\top \mathbf{x} \le u$ into $l\le \mathbf{a}^\top \mathbf{x}$ and $\mathbf{a}^\top \mathbf{x} \le u$ via the \texttt{split-interval} bridge, and then converts each inequality into a \textsf{ScalarAffineFunction}-in-\textsf{EqualTo} constraint and a \textsf{SingleVariable}-in-\textsf{GreaterThan} constraint using the \texttt{slack} bridge, which introduces an additional slack variable $y$. The second reformulation includes an additional step of converting the temporary constraint $\mathbf{a}^\top \mathbf{x} \le u$ into $-u \le -\mathbf{a}^\top \mathbf{x}$ via the \texttt{flip-sign} bridge. Notably, both reformulations add two \textsf{ScalarAffineFunction}-in-\textsf{EqualTo} constraints, two slack variables, and two \textsf{SingleVariable}-in-\textsf{GreaterThan} constraints, but the first reformulation is preferred because it has the least number of transformations.

\subsection{Hyper-graphs and shortest paths}

It is easy to see that as the number of constraint types and bridges increases, the number of different equivalent reformulations also increases, and choosing an appropriate reformulation becomes difficult. 
We overcome the proliferation challenge by posing the question of how to transform a constraint into a set of supported equivalents as a shortest path problem through a directed hyper-graph.

We define our directed hyper-graph $G(N, E)$ by a set of nodes $N$, containing one node $n$ for each possible function-in-set pair, and a set of directed hyper-edges $E$. Each directed hyper-edge $e \in E$, corresponding to a bridge, is comprised of a source node $s(e)\in N$ and a set of target nodes $T(e) \subseteq N$. For each hyper-edge $e$, we define a weight $w(e)$. For simplicity, \texttt{MathOptInterface} chooses $w(e) = 1$ for all $e \in E$, but this need not be the case. In addition, each solver defines a set of \textit{supported} nodes $S$. Finally, for each node $n\in N$, we define a cost function, $C(n)$, which represents the cost of bridging node $n$ into an equivalent set of supported constraints:
$$C(n) = \begin{cases}
0 & n \in S \\
\min\limits_{e\in E\;:\;s(e) = n}\{w(e) + \sum\limits_{n^\prime \in T(e)} C(n^\prime)\}&\text{otherwise}.
\end{cases}$$

In the spirit of dynamic programming, if we can find the minimum cost $C(n)$ for any node $n$, we also obtain a corresponding hyper-edge $e$. This repeats recursively until we reach a terminal node at which $C(n) = 0$, representing a constraint that the solver natively supports. The collection of edges associated with a solution is referred to as a \textit{hyper-path}.

Problems of this form are well studied by \citet{gallo1993directed}, who propose an efficient algorithm for computing $C(n)$ and obtaining the minimum cost edge $e$ associated with each node. Due to the large number of nodes in the hyper-graph, we do not precompute the shortest path for all nodes \textit{a priori}. Instead, we compute $C(n)$ in a \textit{just-in-time} fashion whenever the first constraint of type $n$ is added to the model. Because the computation is performed once per type of constraint, the decision is independent of any constraint data like coefficient values.

The choice of cost function has a significant impact both on the optimal solution and on the computational tractability of the problem.
Indeed, if the cost function is chosen to be the number of different bridges used, the shortest path problem is NP-complete \citep{italiano1989online}.
In the present case, if a bridge is used twice, it makes sense to include its weight twice as well.
This cost function is part of the more general family of \emph{additive cost functions} for which the shortest hyper-path problem can be solved efficiently with a generalization of the Bellman-Ford or Dijkstra algorithms; see \citet[Section~6]{gallo1993directed} for more details.

\subsection{Variable and objective bridges}

In the interest of simplicity, we have described only \textit{constraint bridges}, in which the nodes in the hyper-graph correspond to function-in-set pairs. In practice, there are three types of nodes in $N$:
\textit{constraint} nodes for each pair of function type $f$ and set type $S$ representing $f$-in-$S$ constraints; 
\textit{objective} nodes for each type $f$ representing an objective function of type $f$; and \textit{variable} nodes for each set type $S$ representing variables constrained to $S$. Hyper-edges beginning at a node $n$ can have target nodes of different types.

Objective nodes (and corresponding bridges) allow, for example, conic solvers that support only affine objectives to solve problems modeled with a quadratic objectives by replacing the objective with a slack variable $y$, and then adding a quadratic inequality constraint. If necessary, the quadratic inequality constraint may be further bridged to a second-order cone constraint.

Variable nodes correspond to a concept we call \textit{variables constrained on creation}, and they are needed due to differences in the way solvers handle variable initialization. A na\"ive way of creating variables is to first add $N$ variables to the model, and then add \textsf{SingleVariable} and \textsf{VectorOfVariables} constraints to constrain the domain. This approach works for many solvers, but fails in two common cases: (i) some solvers, e.g., CSDP \citep{borchers1999csdp}, do not support free variables; and (ii) some solvers, e.g., MOSEK \citep{mosek}, have a special type of variable for PSD variables which must be specified at creation time. For example, adding the constraint $X \succeq 0$ to MOSEK after $X$ has been created will result in a bridge that creates a new PSD matrix variable $Y \succeq 0$, and then a set of \textsf{ScalarAffineFunction}-in-\textsf{EqualTo} constraints such that $X = Y$.

Similar to constraints, solvers specify a set of supported variable sets (i.e., so $C(n) = 0$). For most solvers, the supported variable set is the singleton \textsf{Reals}. If the user attempts to add a variable constrained on creation to a set $S$ that is not supported, a bridge first adds a free variable ($x$-in-\textsf{Reals}) and then adds a \textsf{SingleVariable}-in-$S$ constraint. Thus, variable nodes allow solvers such as CSDP to declare that they support only $x$-in-\textsf{Nonnegatives} and not $x$-in-\textsf{Reals}, and they provide an efficient way for users to add PSD variables to MOSEK, bypassing the slack variables and equality constraints that would need to be added if the PSD constraint was added after the variables were created.

It is important to note that constraint and objective bridges are self-contained; they do not return objects that are used in other parts of the model. However, variable bridges \textit{do} return objects that are used in other parts of the model. For example, adding $x\in\textsf{Reals}$ may add two variables $[x^+, x^-] \in \textsf{Nonnegatives}$ and return the expression $x^+ - x^-$ for $x$. The expression $x^+ - x^-$ must then be substituted for $x$ on every occurrence. A detailed description of how this substitution is achieved in code is non-trivial and is outside the scope of this paper.

\subsection{Example}\label{sec:bridge-example}

To demonstrate the combination of the three types of nodes in the hyper-graph, consider bridging a \textsf{ScalarQuadraticFunction} objective function to a solver that supports only:
\begin{itemize}
    \item \textsf{VectorAffineFunction}-in-\textsf{RotatedSecondOrderCone} constraints;
    \item \textsf{ScalarAffineFunction} objective functions; and
    \item Variables in \textsf{Nonnegatives}.
\end{itemize}
As a simple example, we use:
\begin{equation*}
    \begin{array}{rl}
        \min & x^2 + x + 1\\
        \text{s.t.} & x \in \mathbb{R}^1_+.
    \end{array}
\end{equation*}
The first step is to introduce a slack variable $y$ and replace the objective with the \textsf{SingleVariable} function $y$:
\begin{equation*}
    \begin{array}{rl}
        \min & y\\
        \text{s.t.} & x^2 + x + 1 \le y\\
        & x \in \mathbb{R}^1_+\\
        & y\ \mathsf{free}.
    \end{array}
\end{equation*}
However, since the solver supports only \textsf{ScalarAffineFunction} objective functions, the objective function is further bridged to:
\begin{equation*}
    \begin{array}{rl}
        \min & 1y + 0\\
        \text{s.t.} & x^2 + x + 1 \le y\\
        & x \in \mathbb{R}^1_+\\
        & y\ \mathsf{free}.
    \end{array}
\end{equation*}
The second step is to bridge the \textsf{ScalarQuadraticFunction}-in-\textsf{LessThan} constraint into a \textsf{VectorAffineFunction}-in-\textsf{RotatedSecondOrderCone} constraint using the relationship:
\begin{equation*}
    \begin{array}{r c l}
    \frac{1}{2}\mathbf{x}^\top Q \mathbf{x} + \mathbf{a}^\top \mathbf{x} + b \le 0 & \iff & || U\mathbf{x} ||_2^2 \le 2(-\mathbf{a}^\top \mathbf{x} - b)\\
        & \iff & [1, -\mathbf{a}^\top \mathbf{x} - b, U\mathbf{x}] \in \textsf{RotatedSecondOrderCone},
    \end{array}
\end{equation*}
where $Q = U^\top U$. Therefore, we get:
\begin{equation*}
    \begin{array}{rl}
        \min & 1y + 0\\
        \text{s.t.} & [1, -x + y - 1, \sqrt 2x] \in \textsf{RotatedSecondOrderCone}\\
        & x \in \mathbb{R}^1_+\\
        & y\ \mathsf{free}.
    \end{array}
\end{equation*}

Finally, since the solver does not support free variables, a variable bridge is used to convert $y$ into two non-negative variables, resulting in:
\begin{equation*}
    \begin{array}{rl}
        \min & 1y^+ - 1 y^- + 0\\
        \text{s.t.} & [1, -x + y^+ - y^- - 1, \sqrt 2x] \in \textsf{RotatedSecondOrderCone}\\
        & [x, y^+, y^-] \in \mathbb{R}^3_+.
    \end{array}
\end{equation*}
Note how the expression $y^+ - y^-$ is substituted for $y$ throughout the model.

The optimal hyper-path corresponding to this example is given in Figure \ref{fig:hyper-path}. To summarize, the \textsf{ScalarQuadraticFunction} objective node is bridged to a \textsf{SingleVariable} objective node, a $x\in\textsf{Reals}$ variable node, and a $\textsf{ScalarQuadraticFunction}\in\textsf{LessThan}$ constraint node. Then, the \textsf{SingleVariable} objective is further bridged to a \textsf{ScalarAffineFunction} objective node, the $x\in\textsf{Reals}$ variable node is bridged to a $x\in\textsf{Nonnegatives}$ node, and the $\textsf{ScalarQuadraticFunction}\in\textsf{LessThan}$ constraint is bridged to a $\textsf{VectorAffineFunction}\in\textsf{RotatedSecondOrderCone}$ constraint node.

\begin{figure}[!ht]
    \centering
    \resizebox{0.7\textwidth}{!}{%
    \begin{tikzpicture}[align=center]
        \node[draw,dashed, minimum height=2cm] (jump)       at (6, 6) {$\min$\\\textsf{ScalarQuadraticFunction}};
        \node[minimum height=2cm] (lq)         at (6, 3) {$x\in\textsf{Reals}$};
        \node[minimum height=2cm] (con)        at (0, 3) {$\min$\\\textsf{SingleVariable}};
        \node[minimum height=2cm] (nl)         at (12, 3) {\textsf{ScalarQuadraticFunction}\\$\in$\\\textsf{LessThan}};
        \node[draw, minimum height=2cm] (con_solver) at (0, 0) {$\min$\\\textsf{ScalarAffineFunction}};
        \node[draw, minimum height=2cm] (lq_solver)  at (6, 0) {$x\in\textsf{Nonnegatives}$};
        \node[draw, minimum height=2cm] (nl_solver)  at (12, 0) {\textsf{VectorAffineFunction}\\$\in$\\\textsf{RotatedSecondOrderCone}};
        \path[]
        	(jump)	edge [] node [] {} (lq)
        	        edge [] node [] {} (con)
        	        edge [] node [] {} (nl)
        	(lq)	edge [] node [] {} (lq_solver)
        	(con)	edge [] node [] {} (con_solver)
        	(nl)	edge [] node [] {} (nl_solver);
    \end{tikzpicture}}
    \caption{Optimal hyper-path of example in Section \ref{sec:bridge-example}. Dashed box is the objective node we want to add to the model, solid boxes are supported nodes, nodes with no outline are unsupported intermediate nodes, and arcs are bridges.}
    \label{fig:hyper-path}
\end{figure}

\color{revision_color}

\subsection{Benchmarks}\label{sec:bridging_benchmark}

To benchmark the performance of \texttt{MathOptInterface}, and the bridging system in particular, we consider the continuous version of the $P$-median problem used by \citet{hart2011pyomo} to compare Pyomo with AMPL and also by \citet{lubinComputingOperationsResearch2015} when benchmarking an earlier version of JuMP. The model determines the location of $d$ facilities over $N$ possible locations to minimize the cost of serving $M$ customers, where the cost of serving customer $i$ from facility $j$ is given by $c_{ij}$. The decision variable $x_{ij}$ represents the proportion of customer $i$'s demand served by facility $j$, and the decision variable $y_j$ represents the proportion of facility $j$ to open. In practice, this model is usually solved as a mixed-integer linear program with $y_j\in\{0, 1\}$; we consider the continuous relaxation.

We compare two formulations of this problem. The first is a scalar-valued formulation in which we exclusively use functions and sets that are natively supported by the solver GLPK \citep{glpk}:
\begin{equation}\label{eq:scalar_pmedian}
    \begin{array}{rrrll}
        \min & \sum\limits_{i = 1}^M\sum\limits_{j=1}^N c_{ij} x_{ij}\\
        \text{s.t.} & \sum\limits_{j=1}^N x_{ij}& & = 1 & \forall i =1,\ldots,M\\
            & & \sum\limits_{j=1}^N y_{j} & = d \\
            & x_{ij} & -\ y_j & \le 0 & \forall i = 1,\ldots,M,\  j=1,\ldots,N\\
            & x_{ij} & & \ge 0      & \forall i = 1,\ldots,M,\ j = 1,\ldots,N\\
            & & y_j & \in [0, 1]  &  \forall j = 1,\ldots,N.
    \end{array}
\end{equation}
There are $M + 1$ \textsf{ScalarAffineFunction}-in-\textsf{EqualTo} constraints, $M \times N$ \textsf{ScalarAffineFunction}-in-\textsf{LessThan} constraints, $M \times N$ \textsf{SingleVariable}-in-\textsf{GreaterThan} constraints, and $N$ \textsf{SingleVariable}-in-\textsf{Interval} constraints.

The second formulation we consider exclusively uses vector-valued functions that are natively supported by the solver SCS (i.e., in the standard geometric form~\eqref{eq:conic_standard_form}):
\begin{equation}\label{eq:vector_pmedian}
    \begin{array}{rrrl}
        \min & \mathbf{c}^\top \mathbf{x}\\
        \text{s.t.} & A \mathbf{x} & & - \mathbf{1} \in \{0\}^M \\
            & &\mathbf{1}^\top \mathbf{y} & - d  \in \{0\}^1\\
            & B \mathbf{x} &- C \mathbf{y} & + \mathbf{0} \in\mathbb{R}_-^{M\times N}\\
            & I \mathbf{x} & & + \mathbf{0} \in \mathbb{R}_+^{M\times N} \\
            & & I \mathbf{y} & + \mathbf{0} \in \mathbb{R}_+^{N} \\
            & & I \mathbf{y} & - \mathbf{1} \in \mathbb{R}_-^{N},
    \end{array}
\end{equation}
where $\mathbf{c}$, $\mathbf{x}$, and $\mathbf{y}$ are column vectors created by stacking the scalar elements $c_{ij}$, $x_{ij}$, and $y_{j}$, $I$ is the identity matrix, and $A$, $B$, and $C$ are appropriately valued matrices so that model~\eqref{eq:vector_pmedian} is equivalent to model~\eqref{eq:scalar_pmedian}.
There are two \textsf{VectorAffineFunction}-in-\textsf{Zeros} constraints (of dimension $M$ and $1$), two \textsf{VectorAffineFunction}-in-\textsf{Nonpositives} constraints (of dimension $M\times N$ and $N$), and two \textsf{VectorAffineFunction}-in-\textsf{Nonnegatives} constraints (of dimension $M\times N$ and $N$). Note that SCS does not support variable bounds, so we convert variable bounds into \textsf{VectorAffineFunction} equivalents.

Because GLPK natively supports all constraints in the first formulation but none in the second, and SCS supports all constraints in the second but none in the first, these two formulations allow us to test the efficacy of the bridging system. Moreover, this test represents a worst-case bridging scenario in which we have to bridge every constraint.

In addition to the two \texttt{MathOptInterface} models, we also coded a version for GLPK and SCS using their C API directly from Julia, thereby bypassing the overhead of \texttt{MathOptInterface}. Moreover, as an additional baseline, we implemented an equivalent model in CVXPY \citep{diamond2016cvxpy}. To improve the performance of CVXPY, our implementation constructs the $B$ and $C$ matrices directly---at an arguable loss of model readability---because adding the $x_{ij} \le y_j$ constraints individually was many times slower. In total, there are four versions of the problem for both GLPK and SCS, which we refer to as \textit{scalar}, \textit{vector}, \textit{direct}, and \textit{cvxpy}, respectively.

Following \citet{lubinComputingOperationsResearch2015}, in our benchmarks we fixed $d=100$ and $M=100$ and varied $N$. For each problem, we partition the solution time into a \textit{generation} step, in which the data is processed, and a \textit{load} step, in which the data is passed to the solver. While it is natural to compare the overhead of the modeling API with the time to solve the model, i.e., find an optimal solution, the latter depends on algorithmic parameters and tolerances that are challenging to set fairly (the overhead becomes arbitrarily small as the solve time increases). Instead, we consider a worst case by configuring the solvers to terminate immediately (after 1 millisecond for GLPK and 1 iteration for SCS). Hence, we call the time taken to load the data into the solver then immediately return as the \textit{load} time. Our experiments were performed on a Linux system with an Intel Xeon E5-2687 processor, 250 GB RAM, and the following combinations of software: Julia 1.5, Python 3.8, CVXPY 1.1.3, GLPK 4.64, and SCS 2.1.2. The results are shown in Table~\ref{tab:bridging_experiment}. 

\revisionb{Note that CVXPY is a higher-level modeling language than MathOptInterface, and it includes features such as grammar verification for disciplined convex programming that are not present in MathOptInterface. Therefore, there is an understandable and expected performance gap between MathOptInterface and CVXPY. Nevertheless, we include CVXPY as a useful reference point for an open-source modeling language that also reformulates a linear program to the conic form expected by solvers such as SCS, even if CVXPY does more work during the reformulation.}

\begin{table}[!ht]
    \centering
    \begin{tabular}{l r r r r r r r r r}
        \toprule
        & & \multicolumn{4}{c}{GLPK} & \multicolumn{4}{c}{SCS}    \\
        \cmidrule(lr){3-6}\cmidrule(lr){7-10}
        N & &  scalar & vector & direct & cvxpy &  scalar & vector & direct & cvxpy \\
        \midrule
		\multirow{3}{*}{1,000} & generate & 0.05 & 0.14 & 0.05 & 1.51 & 0.48 & 0.05 & 0.07 & 1.38 \\
		& load & 0.19 & 0.46 & 0.06 & 0.26 & 0.71 & 0.18 & 0.11 & 0.10 \\
		& total & 0.24 & 0.59 & 0.11 & 1.77 & 1.19 & 0.23 & 0.18 & 1.48 \\
		\addlinespace
		\multirow{3}{*}{5,000} & generate & 0.23 & 0.72 & 0.31 & 11.49 & 1.85 & 0.11 & 0.54 & 10.89 \\
		& load & 1.40 & 1.84 & 0.41 & 1.55 & 3.99 & 1.51 & 0.82 & 0.72 \\
		& total & 1.63 & 2.57 & 0.72 & 13.04 & 5.84 & 1.62 & 1.36 & 11.61 \\
		\addlinespace
		\multirow{3}{*}{10,000} & generate & 0.69 & 1.35 & 0.65 & 30.27 & 3.98 & 0.28 & 0.91 & 29.31 \\
		& load & 2.49 & 3.64 & 0.70 & 2.77 & 8.28 & 3.17 & 1.81 & 1.42 \\
		& total & 3.19 & 5.00 & 1.34 & 33.04 & 12.26 & 3.45 & 2.72 & 30.73 \\
		\addlinespace
		\multirow{3}{*}{50,000} & generate & 3.32 & 6.58 & 3.32 & 464.19 & 24.65 & 1.17 & 4.55 & 455.42 \\
		& load & 14.04 & 20.90 & 4.09 & 12.33 & 54.45 & 22.38 & 9.39 & 7.26 \\
		& total & 17.36 & 27.48 & 7.41 & 476.52 & 79.09 & 23.55 & 13.94 & 462.68 \\
        \bottomrule
    \end{tabular}
    \caption{$P$-Median benchmark results measuring the overhead of \texttt{MathOptInterface} and bridges with GLPK and SCS for various size $N$. All times in seconds. \textit{scalar} is formulation~\eqref{eq:scalar_pmedian}, \textit{vector} is formulation~\eqref{eq:vector_pmedian}, \textit{direct} is using the C API for each solver, and \textit{cvxpy} an equivalent implementation in CVXPY. \textit{generate} is the time taken to generate the problem data, \textit{load} is the time taken to load the data into the solver and begin the solution process, \textit{total} is the sum.}
    \label{tab:bridging_experiment}
\end{table}

Using \texttt{MathOptInterface} or CVXPY results in overhead compared to using the C API of each solver directly. For CVXPY, this overhead can be a factor of 10 to 20. For \texttt{MathOptInterface}, it depends on whether the bridging system is used. If the bridging system is bypassed (i.e., \textit{scalar} for GLPK and \textit{vector} for SCS), the overhead of \texttt{MathOptInterface} varies between a factor of 1 and 2.5. If the bridging system is used, the overhead is approximately a factor of $4$ for GLPK, and a factor of 4--7 for SCS. Readers may note that the load time for SCS \textit{direct} is always greater than the load time for SCS \textit{cvxpy}. This is because the \textit{direct} mode results in an extra copy in order to shift indices from 1-based Julia arrays to 0-based C arrays. Nevertheless, \textit{direct} always wins on \textit{total} time.

\color{black}

\subsection{Future extensions}

There are many aspects of the shortest path problem that we have not explored in our current implementation. For example, the current implementation in \texttt{MathOptInterface} assigns each bridge a weight of 1.0 in the path. Therefore, the objective of our shortest path problem is to minimize the number of bridges in a transformation. However, it is possible to use other scores for the desirability of a reformulation. For example, in the \textsf{ScalarAffineFunction}-in-\textsf{Interval} example mentioned at the start of this section, the third option may be computationally beneficial since it adds fewer rows to the constraint matrix, even though it adds an extra variable. Therefore, we may assign the corresponding edge in the graph a lower weight (e.g., 0.6).

\revision{As we described in Section \ref{sec:literature},} similar systems for automatically transforming problems have appeared in the literature before, e.g., CVXPY has a similar concept called reductions \citep{cvxpy_rewriting} \revision{and MiniZinc has the ability to redefine constraint transforms \citep{belov2016improved}}. However, \revision{MiniZinc requires the user to manually choose the redefinitions for a particular model, and} since CVXPY targets a fixed set of solvers, it can pre-specify the chain of reductions needed for each solver. On the other hand, our shortest-path formulation enables us to separate the transformation logic from the solvers. The fact that we compute the sequence of transformations \textit{at runtime} additionally makes it possible to use new bridges and sets defined in third-party extensions. 

\revision{Moreover, the bridging system makes it easy to add new solvers to JuMP, because they need only support the minimal set of functions and sets that they natively support. Once this is done, users can use the solver \textit{and} the full modeling power of MathOptInterface. For example, users can solve problems with rotated second order cones and convex quadratic objectives and constraints using a solver that only implements a second order cone.}

\revision{New solvers are already being written to target the flexibility and power of \texttt{MathOptInterface}'s bridging system. One example is \texttt{ConstraintSolver.jl} \citep{constraintsolver}, which extends \texttt{MathOptInterface} by providing constraint programming sets such as \texttt{AllDifferent} and \texttt{NotEqualTo}. A second example is \texttt{Hypatia.jl} \citep{coey2020towards}. \texttt{Hypatia} is a conic interior point algorithm that provides specialized barrier functions for a range of non-standard conic sets (e.g., the relative entropy cone, the log-determinant cone, and the polynomial weighted sum-of-squares cone). If the user forms a model using these new cones and uses \texttt{Hypatia} to solve it, the bridging system is bypassed and the model is passed directly to \texttt{Hypatia}. However, when the user solves the same model with a solver like Mosek \citep{mosek}, then the bridging system reformulates the non-standard cones into an equivalent formulation using only cones that are supported by Mosek. This allows users to focus on modeling their problem in the most natural form, and allows solver authors to experiment with (and easily benchmark and test) novel cones and reformulations.}

\section{A new file format for mathematical optimization}\label{sec:mathoptformat}

As we saw in Section \ref{sec:standard_form}, a model in the MathOptInterface standard form is defined by a list of functions and a list of sets. In this section we utilise that fact to describe a new file format for mathematical optimization problems called \textit{MathOptFormat}. MathOptFormat is a serialization of the MathOptInterface abstract data structure into a JSON file \citep{ecmainternationalJSONDataInterchange2017}, and it has the file-extension \texttt{.mof.json}. A complete definition of the format, including a JSONSchema \citep{json_schema} that can be used to validate MathOptFormat files, is available at \revision{\citet{oscar_dowson_2020_3905406}}. 

In addition, due to the role of file formats in problem interchange, the JSONSchema serves as \textit{the} canonical description of the set of functions $\mathcal{F}$ and sets $\mathcal{S}$ defined in MathOptInterface. We envisage that the schema will be extended over time as more functions and sets are added to the MathOptInterface abstract data structure. 

Importantly, the schema is a concrete representation of the format, and it includes a description of how the sparse vectors and matrices are stored. Moreover, although specified in JSON, this representation utilizes simple underlying data structures such as lists, dictionaries, and strings. Therefore, the format could be ported to a different format (e.g., protocol buffers \citep{protobuf}), \textit{without} changing the basic layout and representation of the data.

\subsection{\revision{Definition and example}}

Rather than rigorously define our format, we shall, in the interest of brevity, explain the main details of MathOptFormat through an example. Therefore, consider the following simple mixed-integer program:
\begin{equation}\label{eq:mof_example}
    \begin{array}{rl}
        \max\limits_{x, y} & x + y\\
        \text{subject to:} & x \in \{0, 1\} \\
                           & y \le 2.
    \end{array}
\end{equation}
This example, encoded in MathOptFormat, is given in Figure \ref{fig:mof_example}.

\begin{figure}[!ht]
    \centering
    \begin{Verbatim}[frame=single,fontsize=\small]
{
"version": {"minor": 5, "major": 0},
"variables": [
  {"name": "x"}, {"name": "y"}
],
"objective": {
  "sense": "max",
  "function": {
    "type": "ScalarAffineFunction",
    "terms": [
      {"coefficient": 1.0, "variable": "x"},
      {"coefficient": 1.0, "variable": "y"}
    ],
    "constant": 0.0
  }
},
"constraints": [{
  "function": {
    "type": "SingleVariable", 
    "variable": "x"
  },
  "set": {"type": "ZeroOne"}
}, {
  "function": {
    "type": "SingleVariable", 
    "variable": "y"
  },
  "set": {
    "type": "LessThan", 
    "upper": 2.0
  }
}] 
}
\end{Verbatim}
    \caption{The complete MathOptFormat file describing model \eqref{eq:mof_example}.}
    \label{fig:mof_example}
\end{figure}

Let us now describe each part of the file in Figure \ref{fig:mof_example} in turn. First, notice that the
file format is a valid JSON file. Inside the document, the model is stored as a single JSON object. JSON objects are key-value mappings enclosed by curly braces (\texttt{\{} and \texttt{\}}). There are four required keys at the top level:
\begin{enumerate}
    \item \textsf{version}: A JSON object describing the minimum version of MathOptFormat needed to parse the file. This is included to safeguard against later revisions. It contains two fields: \textsf{major} and \textsf{minor}. These fields should be interpreted using semantic versioning \citep{semver}. The current version of MathOptFormat is \revision{v0.5}.

    \item \textsf{variables}: A list of JSON objects, with one object for each variable in the model. Each object has a required key \textsf{name} which maps to a unique string for that variable. It is illegal to have two variables with the same name. These names will be used later in the file to refer to each variable.

    \item \texttt{objective}:  A JSON object with one required key:
        \begin{enumerate}
            \item \textsf{sense}: A string which must be \textsf{min}, \textsf{max}, or \textsf{feasibility}.

    If the sense is \textsf{min} or \textsf{max}, a second key \textsf{function}, must be defined:

            \item \textsf{function}: A JSON object that describes the objective function. There are many different types of functions that MathOptFormat recognizes, each of which has a different structure. However, each function has a required key called \revision{\textsf{type}} which is used to describe the type of the function. In this case, the function is \textsf{ScalarAffineFunction}.
        \end{enumerate}
    
    \item \textsf{constraints}: A list of JSON objects, with one element for each constraint in the model. Each object has two required fields:
        \begin{enumerate}
            \item \textsf{function}: A JSON object that describes the function $f_i$ associated with constraint $i$. The function field is identical to the function field in \textsf{objective}; however, in this example, the first constraint function is a \textsf{SingleVariable} function of the variable $x$.
            \item \textsf{set}: A JSON object that describes the set $S_i$ associated with constraint $i$. In this example, the second constraint set is the MathOptFormat set \textsf{LessThan} with the field \textsf{upper}.
        \end{enumerate}
\end{enumerate}

\color{revision_color}

\subsection{Comparison with other formats}

We believe the creation of a new file format is justified because we can now write down problems that cannot be written in any other file format, e.g., programs with exponential cones and complementary constraints. As an illustration, Table \ref{tab:mathoptformat} compares the features supported by MathOptFormat against a number of file formats. The file formats we compare are the Conic Benchmark Format (.cbf) \citep{fribergCBLIB2014Benchmark2016}, the GAMS Scalar Format (.gms) \citep{bussieck2003minlplib}, the LP file format (.lp) \citep{lp_format}, the MPS file format (.mps) \citep{mps_ibm}, the NL file format (.nl) \citep{gayWritingNlFiles2005}, the Optimization Services Instance Language (.osil) \citep{fourerOSiLInstanceLanguage2010}, and the SDPA file format (.sdpa) \citep{sdpa}.

Due to the large number of functions and sets supported by MathOptFormat, we do not compare every combination. Instead we compare a selection of general constraint types for which there are differences between the file formats. Table~\ref{tab:mathoptformat} demonstrates that MathOptFormat generalizes a broad class of problems, from conic formats such as CBF to classical formats for mixed-integer linear programming such as MPS.

\begin{table}[!ht]
    \centering
    \resizebox{\columnwidth}{!}{%
    \begin{tabular}{l c c c c c c c c}
        \toprule
        & .mof.json & .cbf & .gms & .lp & .mps & .nl & .osil & .sdpa\\
        \midrule
Lower and upper bounds on variables & Y & & Y & Y & Y & Y & Y & \\
Integer variables & Y & Y & Y & Y & Y & Y & Y & Y*\\
Binary variables & Y & & Y & Y & Y & Y & Y & \\
Semi-integer and semi-continuous variables & Y & & Y & & Y* & & Y & \\
Linear constraints & Y & Y & Y & Y & Y & Y & Y & Y \\
Quadratic constraints & Y & & Y & Y* & Y* & Y & Y & \\
Second-order cones & Y & Y & Y & & & & Y & \\
Exponential cones & Y & Y & & & & & & \\
Power cones & Y & Y & & & & &  & \\
Positive semidefinite cones & Y & Y & & & & & Y & Y\\
Complementarity constraints & Y & & Y & & & Y & Y & \\
General nonlinear constraints &  & & Y & & & Y & Y & \\
         \bottomrule
    \end{tabular}}
    \caption{Summary of types of constraints supported by various file formats. Y = yes. Y* = some non-standard variations.}
    \label{tab:mathoptformat}
\end{table}

\color{black}

\section{The impact on JuMP and conclusions}\label{sec:jump_1.0}

We created MathOptInterface in order to improve JuMP. Therefore, it is useful to reflect on how JuMP has changed with the introduction of \texttt{MathOptInterface}.


% ...

From a software engineering perspective, the largest change is that 90\% of the code in JuMP was re-written during the transition. In terms of lines of code, \revision{14,341} were added, \revision{10,649} were deleted, \revision{2,994} were modified, and only \revision{2,428} remained unchanged. This represents a substantial investment in engineering time from a large number of individual contributors. In addition, \revision{26,498} lines of code were added to the \texttt{MathOptInterface} package (although, \revision{30}\% of these lines were tests for solvers), and many more thousand lines were added accounting for the more than 20 individual solvers supporting \texttt{MathOptInterface}.

From an architectural perspective, the main change is that instead of representing optimization models using three standard forms, JuMP now represents models using a combination of functions and sets. At the bottom level, instead of solvers implementing one of the three standard forms, they now declare a subset of function-in-set constraint pairs that they natively support, along with supported objective functions and sets for variables constrained on creation. Between these two representations sits the bridging system described in Section \ref{sec:bridges}. Thus, analogous to Figure \ref{fig:jump_architecture}, the JuMP architecture now looks like the diagram in Figure \ref{fig:jump_architecture_moi}.

\begin{figure}[!ht]
    \centering
    \resizebox{0.7\textwidth}{!}{%
    \begin{tikzpicture}[align=center]
        \node[my_node] (jump)       at (4, 6) {JuMP};
        \node[mpb_node, dashed] (mpb) at (4, 3) {};
        \node[] (a1) at (0, 2) {$\min f_1$};
        \node[] (a2) at (2, 2) {$x \in S_2$};
        \node[] (a3) at (4, 2) {$f_3 \in S_3$};
        \node[] (a4) at (6, 2) {$f_4 \in S_4$};
        \node[] (a5) at (8, 2) {$f_5 \in S_5$};
        \node[] (b1) at (0, 4) {$\min f_1$};
        \node[] (b2) at (2, 4) {$x \in S_2$};
        \node[] (b3) at (4, 4) {$f_3 \in S_3$};
        \node[] (b4) at (6, 4) {$f_4 \in S_4$};
        \node[] (b5) at (8, 4) {$f_5 \in S_5$};
        
        \node[my_node] (lq_solver)  at (4, 0) {Solver 2};
        \node[my_node] (con_solver) at (0, 0) {Solver 1};
        \node[my_node] (nl_solver)  at (8, 0) {Solver 3};
        \node[] (mpb_text)  at (-1, 4.65) {\texttt{MathOptInterface}};
        \node[] (mpb_text)  at (4, 3) {... bridges ...};
    
        \path[very thick, <->]
        	(jump)	edge [] node [] {} (b1)
                	edge [] node [] {} (b2)
                	edge [] node [] {} (b3)
                	edge [] node [] {} (b4)
                	edge [] node [] {} (b5)
            (con_solver) edge [] node [] {} (a2)
                    edge [] node [] {} (a1)
            (lq_solver) edge [] node [] {} (a1)
                    edge [] node [] {} (a3)
                    edge [] node [] {} (a4)
            (nl_solver) edge [] node [] {} (a5)
                    edge [] node [] {} (a4);
    \end{tikzpicture}
    }
    \caption{Architecture of the new version of JuMP. \revision{JuMP and the solvers agree on a common set of definitions in the \texttt{MathOptInterface} layer.} JuMP allows users to formulate models using \revision{all} combinations of functions and sets, solvers implement a subset of the complete functionality, and \revision{(if possible)} the bridging system transforms the user-provided model into an equivalent representation of the same model supported by the solver.}
    \label{fig:jump_architecture_moi}
\end{figure}

Despite the major changes at the solver and interface level, little of the user-facing code in JuMP changed (aside from some sensible renaming). An example of a JuMP model using the CPLEX \citep{cplex} optimizer is given in Figure \ref{fig:jump_example}. This deceptively simple example demonstrates many unique features discussed in this paper. The \texttt{t >= 0} variable lower bound is converted into a \textsf{SingleVariable}-in-\textsf{GreaterThan} constraint. The \texttt{Int} tag, informing JuMP that the variable \texttt{t} is an integer variable, is converted into a \textsf{SingleVariable}-in-\textsf{Integer} constraint. The \textsf{SecondOrderCone} constraint is bridged into a \textsf{ScalarQuadraticFunction}-in-\textsf{LessThan} constraint. The \texttt{1 <= sum(x) <= 3} constraint is formulated as a \textsf{ScalarAffineFunction}-in-\textsf{Interval}, and then bridged into a \textsf{ScalarAffineFunction}-in-\textsf{LessThan} constraint and a \textsf{ScalarAffineFunction}-in-\textsf{GreaterThan} constraint. After solving the problem with \texttt{optimize!}, we check that the \texttt{termination\_status} is \texttt{OPTIMAL} before querying the objective value. \revision{Finally, we write out the model to a MathOptFormat file. This model can be loaded in future using \texttt{model = read\_from\_file("example.mof.json")}.}

\begin{figure*}[!ht]
    \centering
    \begin{Verbatim}[fontsize=\small]
using JuMP, CPLEX
model = Model(CPLEX.Optimizer)
@variable(model, t >= 0, Int)
@variable(model, x[1:3] >= 0)
@constraint(model, [t; x] in SecondOrderCone())
@constraint(model, 1 <= sum(x) <= 3)
@objective(model, Min, t)
optimize!(model)
if termination_status(model) == MOI.OPTIMAL
    @show objective_value(model)
end
write_to_file(model, "example.mof.json")
\end{Verbatim}
    \caption{An example \revision{using the version 0.21.3 of JuMP.}}
    \label{fig:jump_example}
\end{figure*}

\subsection{\revision{Other features of \texttt{MathOptInterface}}}\label{sec:other_features}

This paper has described \revision{three main contributions that make writing an algebraic modeling language like JuMP easier: the MathOptInterface abstract data structure; the bridging system for automatically rewriting constraints; and the MathOptFormat file format.} However, the re-write of JuMP \revision{and \texttt{MathOptInterface}} involved many more changes than the ones outlined here. In particular, we have not discussed: 
\begin{itemize}
    \item The API of \texttt{MathOptInterface}, which includes a standardized way to get and set a variety of model and solver attributes \revision{(e.g., names, primal/dual starting points, etc.)}, and the ability to incrementally modify problems in-place (e.g., deleting variables and changing constraint coefficients);
    \item \revision{\texttt{MathOptInterface}}'s \textit{manual} and \textit{automatic} caching modes for solvers that do not support the aforementioned incremental modifications;
    \item JuMP's new \textit{direct} mode, which avoids storing an intermediate copy of the model, bypasses the bridging system, and instead hooks directly into the underlying solver with minimal overhead;
    \item The introduction of a new status reporting mechanism \revision{at the JuMP and \texttt{MathOptInterface} level} featuring three distinct types of solution statuses: \textit{termination status} (Why did the solver stop?), \textit{primal status} (What is the status of the primal solution?), and \textit{dual status} (What is the status of the dual solution?);
    \item \revision{JuMP and \texttt{MathOptInterface}}'s re-vamped support for solver callbacks, offering both \textit{solver-independent} callbacks and \textit{solver-dependent} callbacks, which allow the user to interact with solver-specific functionality; and
    \item \revision{\texttt{MathOptInterface}}'s unified testing infrastructure for solvers, which subjects all solvers to thousands of tests for correctness every time a change is made to the codebase. This testing has revealed bugs and undocumented behavior in a number of solvers.  
\end{itemize}

We leave a description of these changes, and many others, to future work.

For more information on JuMP \revision{and MathOptInterface}, including documentation, examples, tutorials, and source code, readers are directed to \url{https://jump.dev}.

\color{revision_color}
\section*{Supplemental Material}
Supplemental material to this paper is available at \url{https://github.com/jump-dev/MOIPaperBenchmarks}.
\color{black}

\ACKNOWLEDGMENT{
JuMP and MathOptInterface are open-source projects that are made possible by volunteer contributions that include not just writing code, but also finding and reporting bugs, writing and editing documentation, and replying to questions on the community forum. Among the contributors, we owe special thanks to Mathieu Besan{\c c}on, Guilherme Bodin, Chris Coey, Carleton Coffrin, Robin Deits, Twan Koolen, Vitor Nesello, Fran{\c c}ois Pacaud, Robert Schwarz, Issam Tahiri, Mathieu Tanneau, Juan Pablo Vielma, and Ulf Wors{\o}e for their work helping to migrate the JuMP ecosystem to MathOptInterface, which included substantial rewrites of solver interfaces. We thank Ross Anderson for comments on a draft of this paper.

MathOptInterface was conceived during the first JuMP-dev workshop at MIT in June 2017; we thank the MIT Sloan Latin America Office for their sponsorship of this workshop. We thank MIT Sloan, IDEX Bordeaux, and Changhyun Kwon for sponsorship of the 2018 workshop, and we thank NSF under grant OAC-1835443, MISTI under the MIT-Chile UAI
and PUC seed funds, and the Pontificia Universidad Cat\'olica de Chile's Institute of Mathematical and
Computational Engineering and Department of Industrial Engineering for sponsorship of the 2019 workshop.
\revision{B.\ Legat} acknowledges the funding support of an F.R.S.-FNRS fellowship. \revision{O.\ Dowson was supported, in part, by Northwestern University's Center for Optimization \& Statistical Learning (OSL).}
}

\bibliographystyle{informs2014}
\bibliography{main}

\end{document}